\documentclass[a4paper,final]{article} 
\usepackage{amsmath, amsthm, amssymb,dsfont,multicol} 
\usepackage[notcite,notref]{showkeys}  
\usepackage{enumerate}

\newcommand{\pages}[1]{}  

\newenvironment{Proof}[1][]
{\protect\noindent\textbf{Proof#1.} }
{\hfill $\blacksquare$ \\}

\newtheorem{Thm}{Theorem}[section] 

\newtheorem{Exa}[Thm]{Example}
\newtheorem{Lem}[Thm]{Lemma}
\newtheorem{Cor}[Thm]{Corollary}
\newtheorem{Rem}[Thm]{Remark}
\newtheorem{Tra}[Thm]{Transformation}

\newcommand{\tf}[1]{\widetilde{#1}}

\renewcommand{\phi}{\varphi}
\newcommand{\mycom}[1]{}
\newcommand{\dom}{{\mathcal D}}
\newcommand{\One}{\mathds{1}}
\newcommand{\Prob}{\mathbb{P}}

\newcommand{\T}{{\mathcal T}}
\newcommand{\Fskript}{{\mathcal F}}
\newcommand{\Rskript}{{\mathcal R}}
\newcommand{\Askript}{{\mathcal A}}
\newcommand{\Dskript}{{\mathcal D}}

\newcommand{\R}{\mathbb{R}}
\newcommand{\N}{\mathbb{N}}

\newcommand{\E}{\mathbb{E}}
\newcommand{\Bskript}{{\mathcal B}}

\title{Feller Evolution Systems: Generators and Approximation
}
\author{Bj\"orn B\"ottcher\footnote{\texttt{bjoern.boettcher@tu-dresden.de}, TU Dresden, Fachrichtung Mathematik, Institut f\"ur Math. Stochastik, 01062 Dresden, Germany}}
\date{\today}

\begin{document}
\maketitle
\begin{abstract}

A time and space inhomogeneous Markov process is a Feller evolution process, if the corresponding evolution system on the continuous functions vanishing at infinity is strongly continuous. We discuss generators of such systems and show that under mild conditions on the generators a Feller evolution can be approximated by Markov chains with L\'evy increments. 

The result is based on the approximation of the time homogeneous space-time process corresponding to a Feller evolution process. In particular, we show that a $d$-dimensional Feller evolution corresponds to a $d+1$-dimensional Feller process. It is remarkable that, in general, this Feller process has a generator with discontinuous symbol.

\end{abstract}

\bigskip


\bigskip

{\it Secondary:} 
60J25 
35S10 
60J35 

\section{Introduction}

Markov processes with continuous time parameter are used in various applications and often approximations and simulations of these processes are required. In an infinitesimal sense (see the next section) such a process is composed of L\'evy processes. Thus it is natural to try to approximate a Markov process by Markov chains with L\'evy increments. In the case of time homogeneous Markov processes general conditions for such an approximation were given in \cite{BoetSchi2009}. We are going to extend the result to time inhomogeneous Markov processes. The main tool is the transformation of a time inhomogeneous Markov process to a time homogeneous Markov process, which will be recalled in Section \ref{sectransfomation}. The transformed process is usually called the corresponding \textit{space-time process} and the transformation was already used by Doob \cite[p.\ 226]{Doob1955} and Dynkin \cite[Section 4.6\pages{, p.\ 102}]{Dynk1961}. We will derive a necessary and sufficient condition for the transformed process to be a Feller process. Furthermore, the generator of this Feller process is analyzed. In Section \ref{secapprox} the approximation is presented and discussed.

\bigskip

For better readability column vectors will be written as rows and a process $(X_t)_{t\in \R},$ an evolution system $(U(s,t))_{s,t\in \R, s\leq t}$ and a semigroup $(T(t))_{t\in\R}$ will be often denoted just by $X_t,$ $U(s,t)$ and $T(t),$ respectively. The Borel measurable functions on $\R^d$ will be denoted by $B(\R^d)$ and the continuous functions by $C(\R^d)$. The subscripts $c$, $b$ and $\infty$ denote functions with compact support, bounded functions and functions vanishing at $\infty$, respectively; furthermore, a superscript indicates the number of existing derivatives. 
The uniform norm is denoted by $\|.\|_\infty.$

\section{Markov processes and Generators}
Let $X_t$ be an $\R^d$ valued time (and space) inhomogeneous Markov process on the probability space $(\Omega, \Askript, \Prob).$ Then the corresponding evolution system 
\[U(s,t) f(x) := \E(f(X_t)|X_s=x),\ s\leq t,\ s,t\in\R\]
is well defined on $B_b(\R^n).$ The linear operators $U(s,t)$ are positivity preserving and satisfy $U(s,t)1 = 1$, $U(s,s)=\mathrm{id}$ and the evolution property $U(s,t) = U(s,r)U(r,t)$ for $s\leq r \leq t.$

Such families of operators are well studied in the literature, e.g.\ Yosida \cite[Section XIV.4\pages{, p.\ 430}]{Yosi1971}, Pazy \cite[Chapter 5\pages{, p.\ 127}]{Pazy1983}. The following definitions are analogous to Gulisashvili and van Casteren \cite[Section 2.3\pages{, p.\ 106}]{GuliCast2006}, who use the term \textit{backward propagator} for an evolution system.

Corresponding to an evolution system a family of right generators is given by  
\begin{equation}\label{generator}
A^+_s f :=  \lim_{h\downarrow 0} \frac{U(s,s+h)f-f}{h} \ \ \ \textnormal{ for each }s\in \R
\end{equation}
which is defined for all $f\in C_\infty(\R^d)$ such that the limit exists in a strong sense (i.e.\ with respect to $\|\cdot\|_\infty$). In this case we write $f\in \dom(A_s^+).$ If one weakens \eqref{generator} to a pointwise limit, the corresponding operator is called \emph{extended pointwise generator} (a notion which will be of importance in Theorem \ref{calgenerator}). Analogously the left generators are defined by
\begin{equation*}
A_s^- f :=  \lim_{h\downarrow 0} \frac{U(s-h,s)f-f}{h} 
\end{equation*}
on $\dom(A^-_s).$

The family of operators $U(s,t)$ is strongly continuous, if for each $v, w \in \R , v\leq w$
\begin{equation}\label{strongcont}
\lim_{\substack{(s,t)\to (v,w)\\ s\leq t}} \|U(s,t)f - U(v,w)f\|_\infty = 0.
\end{equation}
Note that a family of linear operators $U(s,t)$ on $C_\infty$ satisfying \eqref{strongcont}, $\|U(s,t)f \|_\infty \leq \|f\|_\infty$, $U(s,t)f(x) = U(s,r)U(r,t)f(x)$ for $s\leq r \leq t$ and $U(s,t) f \geq 0$ for $f \geq 0$ is called a \textit{Feller evolution system}.

We denote by $\frac{d}{dt}^+$ ($\frac{d}{dt}^-$) the right (left) derivative. Thus the evolution property leads to the following evolution equations corresponding to the process:
\begin{eqnarray}
\label{forward}
\frac{d}{dt}^\pm U(s,t) & = &U(s,t)A^\pm_t, \\
\label{backward}
\frac{d}{ds}^\pm U(s,t) & = &-A^\pm_s U(s,t).
\end{eqnarray}
Equation \eqref{forward}$^+$ is called forward equation and \eqref{backward}$^-$ is called backward equation. Note that only in the case of the backward equation it makes sense to talk about solutions of the corresponding initial value problem. In the case of the forward equation one can only consider fundamental solutions due to the interchanged order of $U$ and $A.$ These equations are equivalent to the Kolmogorov equations if the corresponding process has transition densities $p(s,x;t,y).$ Then
$$U(s,t)f(x) = \int f(y) p(s,x;t,y)\ dy =: \langle f,p(s,x;t,.)\rangle_{L^2}$$
holds and thus the forward equation reads as
$$\big\langle f,\tfrac{d}{dt}p(s,x;t,.)\big\rangle_{L^2} = \big\langle f,{A_t^{+}}^\star p(s,x;t,.)\big\rangle_{L^2},$$
where ${A_t^{+}}^\star$ is the (formal) adjoint of $A_t^{+}$.

The operators $A_s^+$ (resp. $A_s^-$) satisfy the {\em positive maximum} principle, i.e.\ for $f\in \dom(A^+_s)$ the following implication holds: If there exists $x_0 \in \R^d$ with $f(x_0) = \sup_{x\in\R^n} f(x) \geq 0,$ then
\[ A^+_sf(x)\big|_{x=x_0} \leq 0.\]

This property of $A_s^+$ (resp. $A_s^-$) is a consequence of \eqref{generator} and the fact, that for $f\in\dom(A_s^+)$ (resp. $A_s^-$) attaining its positive maximum at some point $x_0$ the following inequality holds:
\[U(s,t) f(x_0) \leq U(s,t)f^+(x_0) \leq \|f^+\|_\infty = f(x_0), \textnormal{ where }f^+:= f \One_{\{f\geq 0\}}.\]

Therefore, if the set $C_c^\infty(\R^d)$ is a subset of $\dom(A_s^+)$ (resp. $\dom(A_s^-)$), we know by Courr\`{e}ge \cite{Cour66} (see Jacob \cite[Section 4.5]{Jaco2001}) that $-A_s^+$ (resp. $-A_s^-$) on $C_c^\infty$ is a pseudo-differential operator with {\em continuous negative definite symbol}, i.e.\ it admits the representation
\begin{equation}
\label{eq-symbol}
-A_s^+ f(x) = (2\pi)^{-\frac{d}{2}} \int_{\R^d} e^{ix\xi} q_+(s,x,\xi) \hat f(\xi) \ d\xi,
\end{equation}
where $\hat f(\xi) = (2\pi)^{-\frac{d}{2}} \int_{\R^d} e^{-ix\xi} f(x) \ dx$ denotes the Fourier transform of $f$ and $q_+(s,x,\cdot)$ is for fixed $(s,x)$ a continuous negative definite function in the sense of Berg and Forst \cite{BergFors75}. Conversely, for any operator defined via \eqref{eq-symbol} the function $q_+$ is called the \emph{symbol} of the operator.

An explicit construction of Feller evolutions for a given symbol can be found in \cite{Boet2008}. Another option is to construct a Feller process with a constant drift coordinate and consider the process of the remaining coordinates as Feller evolution (see Theorem \ref{feller} and Corollary \ref{approxconstruct} below). For a survey of constructions of Feller processes for a given symbol see \cite{JacoSchi2001}. 

Note that for a time homogeneous evolution system (i.e.\ $U$ such that $U(s,s+h)=U(s-h,s)$ for all $s,h\geq 0$) the left and right generators coincide and do not depend on time. Thus $T(h):=U(s,s+h)$ defines a Feller semigroup with generator $A:=A^-_s=A^+_s.$

In general, however, the left and the right generator do not coincide as the following example illustrates.

\begin{Exa} \label{exa-lrgen}
Consider a process on $\R$ which drifts with slope $\alpha>0$ until time $s_0$ and afterwards with slope $\beta>0$, $\alpha\neq\beta.$ Thus the process started in $x\in\R$ at time $t_0$ is given by
$$X_t = \begin{cases}x + \alpha\cdot (t-t_0) &, t<s_0 \\ x+ \alpha\cdot (s_0-t_0) + \beta\cdot (t-s_0) &, t\geq s_0\end{cases} \quad \text{ for }t\geq t_0,$$
and for $f\in C_c^\infty(\R)$
$$A_{s_0}^-f(x) = \lim_{h\downarrow 0} \frac{\E(f(X_{s_0})|X_{s_0-h}=x)-f(x)}{h} = \lim_{h\downarrow 0} \frac{f(x+\alpha h) - f(x)}{h} =  \alpha f'(x),$$
$$A_{s_0}^+f(x) = \lim_{h\downarrow 0} \frac{\E(f(X_{s_0+h})|X_{s_0}=x)-f(x)}{h} = \lim_{h\downarrow 0} \frac{f(x+\beta h) - f(x)}{h} =  \beta f'(x).$$
Moreover the symbol corresponding to $A_s^+$ as in \eqref{eq-symbol} is given by
$$q_+(s,x,\xi) = - i l(s)\xi \text{ with } l(s) := \begin{cases} \alpha &, s<s_0,\\ \beta &, s\geq s_0. \end{cases}$$
(See Remark \ref{exa-lrgen-spacetime} for the space-time transformation of this process.)
\end{Exa}

The following lemma gives some condition for the left and right generator to coincide.

\begin{Lem} \label{gen-coincide}
Fix $s$ and select $f\in \Dskript(A_s^+)$ such that there exists some $\delta>0$ with $f\in \bigcap_{r\in (s-\delta,s]}\Dskript(A_r^+).$ If $A_r^+$ exists uniformly for $r\in (s-\delta,s],$ i.e.
\begin{equation} \label{unifgen}
\lim_{h\downarrow 0} \sup_{r\in(s-\delta,s]} \left\|\frac{U(r,r + h)f - f}{h} - A^+_r f\right\|_\infty = 0
\end{equation}
and $r \mapsto A_r^+$ is strongly continuous from the left in $s$, i.e.
\begin{equation}
\label{eq-leftcontgen}
\lim_{h \downarrow 0} \|A_{s-h}^+ f - A_s^+ f\|_\infty = 0,
\end{equation}
then $f\in \Dskript(A_s^-)$ and 
$$A_s^+ f=A_s^- f.$$
\end{Lem}
\begin{Proof}
We have for $h<\delta$
\begin{equation*}
\begin{split}
\lim_{h\downarrow 0} & \left\|\frac{U(s-h,s) f - f}{h} - A_s^+f\right\|_\infty\\
& \leq \lim_{h\downarrow 0} \left\|\frac{U(s-h,s-h+h) f - f}{h} - A_{s-h}^+ f\right\|_\infty + \lim_{h \downarrow 0} \|A_{s-h}^+ f - A_s^+ f\|_\infty\\
& \leq \lim_{h\downarrow 0} \sup_{r\in (s-\delta,s]} \left\|\frac{U(r,r + h)f - f}{h} - A^+_r f\right\|_\infty + 0\\
& = 0
\end{split}
\end{equation*}
and the result follows. 
\end{Proof}

The assumptions of Lemma \ref{gen-coincide} are, for example, satisfied if the generator is a pseudo-differential operator whose symbol continuously depends on time (implying \eqref{eq-leftcontgen}) and has bounded coefficients in the sense of \eqref{fellerevobound} below (implying \eqref{unifgen}).

\section{Transformation of time inhomogeneous\\ Markov processes}
\label{sectransfomation}

To transform an $\R^d$ valued time inhomogeneous Markov process $X_t$ defined on $(\Omega,\Askript,\Prob)$ into a time homogeneous Markov process $\tf{X}_t$ defined on $(\tf{\Omega},\tf{\Askript},\tf{\Prob})$ we follow \cite[Section 8.5.5]{Went79}. For generality in this section $\T$ will denote the time set on which $X_t$ is defined, i.e.\ so far we considered $\T = \R$. But $\T = [0,\infty)$ would also be possible. Note that in both cases the transformed process $\tf{X}_t$ will always be defined only on the time set $[0,\infty).$

For $X_t$ there exists a transition function $P:\T \times \R^d \times\T  \times \Bskript \to [0,1]$ such that for each $s,t\in \T , s\leq t, x\in \R^d, B\in \Bskript$ the function $P(s,\cdot;t,B)$ is measurable, $P(s,x;t,\cdot)$ is a probability measure, $P(s,x;s,B) = \One_B(x)$ and $P(s,X_s;t,B)= \Prob(X_t\in B|X_s).$ Furthermore, since $X_t$ is a Markov process also the Chapman Kolmogorov equations 
$P(s,x;t,B) = \int_{\R^d} P(r,y;t,B)\ P(s,x;r,dy)$ hold for $s\leq r \leq t$ and $x,y \in \R^d.$

The standard way to define the transformed process is:

\begin{Tra} Let $X_t$ be as above.
\begin{itemize}
\item {\bf New state space}: $\T  \times \R^d$ with elements $\tf{x} := (s,x),$ $s\in\T ,$ $x\in\R^d.$\\ On this space we consider the $\sigma$-algebra $\tf{\Bskript}$ consisting of all sets $B \subset \T \times\R^d$ such that for all $s\in\T $ the cuts $B_s:=\{x \,:\, (s,x) \in B\}$ are elements of the Borel $\sigma$-algebra on $\R^d$.
\item {\bf New sample space}: $\tf{\Omega}:= \T  \times \Omega$ with elements $\tf{\omega} := (s,\omega),$ $s\in\T ,$ $\omega \in \Omega,$ and the $\sigma$-algebra $\tf{\Askript}:= \{A \subset \tf{\Omega} \,:\, A_s\in \Askript, \forall s \in \T \}$ where $A_s := \{\omega \,:\, (s,\omega) \in A\}.$
\item {\bf Space-time process}: 
$$\tf{X}_t(\tf{\omega}) = \tf{X}_t(s,\omega) := \left(s+t, X_{s+t}(\omega) \right), \ t\in[0,\infty)$$
with the probability measure defined for $A \in \tf{\Askript}$ and $\tf{x} \in \T  \times \R^d$ by
$$\tf{\Prob}_{\tf{x}}(A) =  \tf{\Prob}(A|\tf{X}_0 = (s,x)) := \Prob(A_s |X_s = x),$$
i.e.\ the transition probabilities are given by 
$$\tf{\Prob}(\tf{X}_t\in B | \tf{X}_0 = \tf{x}) = \tf{\Prob}(\tf{X}_t\in B | \tf{X}_0 = (s,x) ) = \Prob( X_{s+t} \in B_{s+t} | X_s = x)$$
where $B \in \tf{\Bskript}, \tf{x} \in \T  \times \R^d$, and thus the transition function is defined by
$$\tf{P}(t,\tf{x},B):= P(s,x;s+t,B_{s+t}).$$
\end{itemize} 
\end{Tra}

In the transformation the change of the probability space might seem counterintuitive, since the process is extended by adding a deterministic drift in a further dimension but no further randomness is introduced. Nevertheless it is canonical, if one recalls the construction of Markov processes using Kolmogorov's theorem.

To see that the new process is a Markov process denote by $\Fskript_t$ its filtration and note that for each $t\in[0,\infty)$, $\tf{x}\in\T\times\R^d$, $B\in\tf{\Bskript}$ the function $\tf{x} \mapsto \tf{P}(t,\tf{x},B)$ is measurable, since for $\tf{x} = (s,x)$ the function is measurable in $x$ and the $\sigma$-Algebra corresponding to $s$ is the power set of $\T$. Furthermore, $\tf{P}(t,\tf{x},.)$ is a probability measure, $\tf{P}(0,\tf{x},B) = \One_B(\tf{x})$ holds and for all $t,r\in\T$ and $\tf{x},\tf{y}\in \T\times \R^d$ the Chapman Kolmogorov equation 
\begin{equation*}
\begin{split}
\int \tf{P}(r,\tf{y},B) \ \tf{P}(t,\tf{x},d\tf{y}) &= \int_{\R^d} P(s+t,y; s+t+r, B_{s+t+r})\ P(s,x;s+t,\ dy)\\
& = P(s,x;s+t+r,B_{s+t+r})\\
& = \tf{P}(t+r,\tf{x}, B)
\end{split}
\end{equation*}
holds.
Hence the process is a Markov process if and only if 
$$\tf{\Prob}_{\tf{x}}(\tf{X}_{t+h}\in B | \Fskript_t) = \tf{P}(h,\tf{X}_t,B)$$
i.e.\ for all $A\in \Fskript_t$
$$\tf{\Prob}_{\tf{x}}(A \cap \{\tf{X}_{t+h} \in B\}) = \int_A \tf{P}(h,\tf{X}_t(\tf{\omega}), B)\ \tf{\Prob}_{\tf{x}}(d\tf{\omega}).$$
This equality holds since the transformation given above and the Markov property of $X_t$ yield
\begin{equation*}
\begin{split}
\int_A \tf{P}(h,\tf{X}_t(\tf{\omega}), B)\ & \tf{\Prob}_{\tf{x}}(d\tf{\omega}) \\
& = \int_{A_s} P(s+t,X_{s+t}(\omega); s+t+h, B_{s+t+h})\ \Prob(d\omega|X_s=x)\\
& = \Prob(A_s \cap \{X_{s+t+h} \in B_{s+t+h}\} | X_s =x )\\
& = \tf{\Prob}_{\tf{x}}(A \cap \{\tf{X}_{t+h} \in B\}).
\end{split}
\end{equation*}

Thus $\tf{X}_t$ is a Markov process and therefore there exists a corresponding semigroup $T(t)$ on $f\in B_b(\T  \times \R^d)$ given by
\begin{equation} \label{semigroup}
 T(t)(\tf{x}) = \tf{\E}(f(\tf{X}_t)|\tf{X}_0 = \tf{x}).
\end{equation}

The following theorem provides a necessary and sufficient condition for $T(t)$ to be a Feller semigroup, i.e.\ a strongly continuous positivity preserving contraction semigroup on $C_\infty$. 
\begin{Thm} \label{feller}
Let $X_t$ be a Markov process with corresponding evolution system $U(s,t).$ Furthermore, let $\tf{X}_t$ be the time homogeneous transformation (as defined above) of $X_t,$ and $T(t)$ be the semigroup associated with $\tf{X}_t$ as in \eqref{semigroup}. 
Then the following statements are equivalent:
\begin{enumerate}[i)]
\item $(U(s,t))_{s,t\in \T, s\leq t}$ is a Feller evolution system on $C_\infty(\R^d)$,
\item $(T(t))_{t\geq 0}$ is a Feller semigroup on $C_\infty(\T \times \R^d)$.
\end{enumerate}
\end{Thm}
\begin{Proof}
Note that 
$$C_\infty(\T  \times \R^d) = \{ f\in C(\T  \times \R^d) \,:\, \lim_{|\tf{x}|\to \infty} f(\tf{x}) = 0\}.$$
Let $f\in C_\infty(\T  \times \R^d)$ and define $g_s(x) := f(s,x)= f(\tf{x})$ for all $\tf{x}=(s,x) \in \T  \times \R^d$ then the semigroup has the representation
\begin{equation}\label{translation}
\begin{split}
T(t)f(\tf{x}) &= \tf{\E}(f(\tf{X_t})| \tf{X_0} = \tf{x} ) \\
&= \tf{\E}( f(\tf{X_t})| \tf{X_0} = (s,x)) \\
&= \E(f(s+t,X_{s+t})| X_s = x) \\
&= \E(g_{s+t}(X_{s+t})| X_s = x)\\
&= U(s,s+t) g_{s+t}(x).
\end{split}
\end{equation}

\noindent i)$\Rightarrow$ii): Clearly $T(t)$ is a positivity preserving contraction semigroup on $B_b$. Thus it remains to show
that 
\begin{enumerate}
\item $T(t)$ maps  $C_\infty(\T  \times \R^d)$ into $C_\infty(\T  \times \R^d)$, 
\item $T(t)$ is strongly continuous on $C_\infty(\T  \times \R^d)$.
\end{enumerate}
First note, that since $U(s,t)$ is strongly continuous, it is also locally uniform strongly continuous, i.e.\ for each compact $K\subset \T $
\begin{equation}\label{lusc}
\lim_{\substack{(s,t)\to (v,w)\\ s\leq t}} \|U(s,t)g - U(v,w)g\|_\infty = 0
\ \ \text{ uniformly for } v,w \in K,\ v\leq w.
\end{equation}

\noindent 
The first step of proving 1. is to show that $\tf{x} \mapsto T(t)f(\tf{x})$ is continuous.
Let $\tf{x}=(s,x), \tf{y}=(r,y)\in \T \times \R^d$ with $\tf{y}$ fixed. Then
\begin{equation*}
\begin{split}
|T(t)f(\tf{x}) - T(t)f(\tf{y})| & = |U(s,s+t) g_{s+t}(x) - U(r,r+t) g_{r+t}(y)|\\
& \leq |U(s,s+t) (g_{s+t} - g_{r+t})(x))| \\
&\phantom{ \leq |}+ |U(s,s+t)g_{r+t}(x) - U(r,r+t)g_{r+t}(x)|\\
&\phantom{ \leq |}+ |U(r,r+t) g_{r+t}(x) - U(r,r+t) g_{r+t}(y)|
\end{split}  
\end{equation*}
holds, where the first term can be estimated by $\|g_{s+t}(.)-g_{r+t}(.)\|_\infty$ using the contraction property of $U(s,t)$ and the second term converges by the local uniform strong continuity \eqref{lusc}. Thus each of these terms is smaller than $\frac{\varepsilon}{3}$ for $|s-r|<\delta_1$ for some $\delta_1>0$. Furthermore, $r$ and $t$ being fixed the function $x \mapsto U(r,r+t)g_{r+t}(x)$ is continuous. Hence also the last term gets smaller than $\frac{\varepsilon}{3}$ for $|x-y|<\delta_2$ for some $\delta_2>0$. Thus taking $|\tf{x} - \tf{y}|<\delta_1 \land \delta_2$ yields the continuity.

The next step is to show that $T(t)f(\tf{x}) \xrightarrow{|\tf{x}|\to \infty} 0.$ Let $\varepsilon >0.$ It holds that
\begin{equation}\label{large-s-vanish}
|T(t)f(\tf{x})| = |U(s,s+t) g_{s+t}(x)| \leq \sup_x |g_{s+t}(x)| = \sup_x |f(s+t,x)|,
\end{equation}
and note that $|\tf{x}|^2=|s|^2+|x|^2,$ i.e.\ for $|\tf{x}| \to \infty$ at least one of $|s|$ and $|x|$ is large. 

Since $f\in C_\infty(\T  \times \R^d)$ there exists $R_1=R_1(t)$ such that $|f(s+t,x)|< \varepsilon$ uniformly in $x$ for $|s|^2 > \frac{R_1}{2}.$

Otherwise, if $|s|^2 \leq \frac{R_1}{2}$, let $h_t(x):=\sup_{|s|^2 \leq {\frac{R_1}{2}}} |g_{s+t}(x)|$ and note that $h_t \in C_\infty(\R^d).$ Thus
$$|T(t)f(\tf{x})| = |U(s,s+t) g_{s+t}(x)| \leq U(s,s+t)h_t(x).$$
By the uniformity of \eqref{lusc} and the Heine-Borel theorem the set $\{s\in\T : |s|^2\leq \frac{R_1}{2}\}$ can be covered by equally sized balls with centres in some finite set $\Rskript_1 \subset \{s\in\T : |s|^2\leq \frac{R_1}{2}\} $ such that
$$\left\| \sup_{|s|^2\leq \frac{R_1}{2}}U(s,s+t)h_t - \max_{r\in\Rskript_1} U(r,r+t) h_t\right\|_\infty \leq \frac{\varepsilon}{2}.$$
Since $\Rskript_1$ is finite and $U(s,t)$ is an evolution system on $C_\infty$ there exists an $R_2$ such that for $|x|^2>\frac{R_2}{2}$ 
$$\max_{r\in \Rskript_1} U(r,r+t) h_t(x) < \frac{\varepsilon}{2}.$$
Hence for $|\tf{x}| > R:= \sqrt{R_1 \lor R_2}$ either $|s|^2 > \frac{R_1}{2}$ and thus \eqref{large-s-vanish} implies the result or $|s|^2\leq \frac{R_1}{2}$ and $|x|^2 > \frac{R_2}{2}$ and therefore
$$|T(t)f(\tf{x})| \leq \sup_{|s|^2\leq\frac{R_1}{2}} U(s,s+t)h_t(x) < \varepsilon.$$

\noindent To show 2. note that 
\begin{equation*} 
\begin{split}
|T(t) f( &\tf{x}) - f(\tf{x}) | = |U(s,s+t) g_{s+t}(x) - g_{s}(x)|\\
&\leq |U(s,s+t)g_{s+t}(x) - U(s,s+t) g_{s}(x)| + |U(s,s+t) g_s(x) - g_s(x)|\\
&\leq \|g_{s+t} - g_s\|_{\infty}  + |U(s,s+t) g_s(x) - g_s(x)|\\
&\leq \left[\sup_{s,x} \big|f(s+t,x) - f(s,x)\big|\right] + |U(s,s+t) g_s(x) - g_s(x)|.
\end{split}
\end{equation*}
holds, where the first term converges to 0 as $t\to 0$ due to the uniform continuity of $f.$ For the second term fix $\varepsilon >0$. Then there exists (analogous to \eqref{large-s-vanish} and the reasoning thereafter) an $R>0$ such that  
$$\sup_x \sup_{|s|>R} |U(s,s+t)g_s(x) - g_s(x)|< \frac{\varepsilon}{3}$$
holds uniformly for all $t<1$. Furthermore, $r\mapsto U(s,s+t)g_r(x) - g_r(x)$ is equicontinuous since 
$$|U(s,t+s) g_r(x) -g_r(x) - U(s,t+s)g_q(x) + g_q(x)| \leq 2 \|g_r - g_q\|_\infty$$
and therefore we find, as above, a finite set $\Rskript \subset \{r\in\T : |r|\leq R\}$ such that for all $x\in\R^d$
$$\sup_{|r|\leq R} |U(s,s+t) g_r - g_r(x)| \leq \max_{r\in \Rskript} |U(s,s+t) g_r(x) - g_r(x)| + \frac{\varepsilon}{3}.$$
Since $\Rskript$ is finite and $U(s,t)$ satisfies \eqref{lusc} there exists a $\delta>0$ such that
$$\sup_{|s|\leq R} \max_{r\in \Rskript} |U(s,s+t) g_r(x) - g_r(x)| < \frac{\varepsilon}{3}\ \ \ \text{ for } t<\delta.$$

Putting the above together yields for $t < \delta \land 1$
\begin{equation*}
\begin{split}
\sup_x \sup_s |U(s,s+t) g_s(x) - g_s(x)| & \leq \sup_x \sup_{|s|\leq R} |U(s,s+t) g_s(x) - g_s(x)| + \frac{\varepsilon}{3}\\
& \leq \sup_{|s| \leq R} \max_{r\in \Rskript} |U(s,s+t) g_r(x) - g_r(x)| + \frac{2 \varepsilon}{3}\\
& < \varepsilon.
\end{split}
\end{equation*}
Thus $T(t)$ is strongly continuous.

\bigskip

\noindent ii)$\Rightarrow$i): $X_t$ is a Markov process, hence $U(s,t)$ is a positivity preserving contraction evolution system on $B_b$. Equation \eqref{translation} with $f(\tf{x}) : = g(x),$ $g\in C_\infty(\R^d)$ reads as
$$T(t)f(\tf{x}) = U(s,s+t)g(x)$$
and thus $U(s,s+t): C_\infty(\R^d) \to C_\infty(\R^d).$ 

\newcommand{\vect}[2]{(#1,#2)}

The function $\tf{x} \mapsto T(t)f(\tf{x})$ is uniformly continuous and thus $s\mapsto T(t)f\vect{s}{x}$ is equicontinuous (w.r.t.\ $x\in\R^d$). Finally \begin{equation*}
\begin{split}
|U(s,t)g(x)-U(v,w)g(x)| & = \left|T(t-s)f\vect{s}{x} - T(w-v)f\vect{v}{x}\right|\\
& \leq \left|T(t-s)f\vect{s}{x} - T(w-v)f\vect{s}{x}\right| \\
& \phantom{\leq} \ + \left|T(w-v)f\vect{s}{x} - T(w-v)f\vect{v}{x}\right|\\
& \leq \sup_{\tf{x}\in \T\times \R^d} |T(t-s) f(\tf{x}) - T(w-v)f(\tf{x})|\\
& \phantom{\leq} \ + \sup_{x \in\R^d}\left|T(w-v)f\vect{s}{x} - T(w-v)f\vect{v}{x}\right|
\end{split}
\end{equation*} 
yields the strong continuity.
\end{Proof}

Now it is straightforward to calculate the generator of $T(t)$: 

\begin{Thm} \label{calgenerator} 
Let $X_t$ be a Feller evolution with evolution system $U(s,t)$ and right generators $A^+_s.$ Furthermore, let $\tf{X}_t$ be its time homogeneous transformation with associated semigroup $T(t)$ as in \eqref{semigroup}. Then the (extended pointwise) generator $L$ of $T(t)$ is given for all $f\in C_\infty(\T  \times \R^d)$ satisfying 
\begin{itemize}
\item $f(.,x) \in C^1(\T ) \text{ for all } x\in \R^d,$
\item $f(s,.) \in \Dskript(A^+_s) \text{ for all } s\in \T $
\end{itemize}
by
\begin{equation}\label{fellergen}
Lf(\tf{x}) = \frac{\partial}{\partial s} f(s,x) + A_s^+ g_s(x) \ \ \text{ where } \tf{x} = (s,x) \text{ and } g_s(x) = f(s,x).
\end{equation}
\end{Thm}
\begin{Rem}
Note that \eqref{fellergen} does not imply that the given $f$ is in the domain of the generator $L$. For this one would have to ensure that $Lf$ is in $C_\infty(\T \times \R^d)$. See Remark \ref{exa-lrgen-spacetime} for further discussion and Lemma \ref{lem-dom} for a sufficient condition.
\end{Rem}
\begin{Proof}[ of Theorem \ref{calgenerator}]
Using \eqref{translation} yields 
\begin{equation*}
\begin{split}
\frac{T(t) f(\tf{x}) - f(\tf{x})}{t} = &\frac{\E(f(s+t,X_{s+t}) | X_s = x) - \E(f(s,X_{s+t})|X_s=x)}{t} \\
&+ \frac{U(s,s+t)g_s(x) - g_s(x)}{t}
\end{split}
\end{equation*}
where the second term converges for $t\to 0$ to $A^+g_s(x).$

For the first term set $f^{(1,0)}(s,x):= \frac{\partial}{\partial s}f(s,x)$ and note that 
\begin{eqnarray*}
\E\bigg( \frac{f(s+t,X_{s+t})-f(s,X_{s+t})}{t} \bigg| X_s=x\bigg) \\
= \E\left(\frac{1}{t} \int_s^{s+t} f^{(1,0)}(r,X_{s+t})\,dr\bigg|X_s=x\right).
\end{eqnarray*}
Furthermore, 
\begin{eqnarray*}
\bigg|\E\left(\frac{1}{t} \int_s^{s+t} f^{(1,0)}(r,X_{s+t})\,dr\bigg|X_s=x\right) - \E(f^{(1,0)}(s,X_{s+t})|X_s=x)\bigg| \\
\leq \sup_y \sup_{r \in[s,s+t]}|f(r,y) - f(s,y)|
\end{eqnarray*}
vanishes for $t \to 0$ since $f$ is uniformly continuous. Finally, defining the function $h_s(x) := f^{(1,0)}(s,x) $ and using strong continuity yields
$$\E(f^{(1,0)}(s,X_{s+t})|X_s=x) = U(s,s+t)h_s(x) \xrightarrow{t\to 0} h_s(x) = \frac{\partial}{\partial s} f(s,x),$$
which proves the statement. 
\end{Proof}

Furthermore if $C_c^\infty$ is in the domain of all right generators, then the operator $L$ has a representation as pseudo-differential operator:  
\begin{Cor} \label{calsymbol} Let $\T = \R$.
If $A^+_s$ is defined on the test functions $C_c^\infty(\R^d)$ and thus can be represented (cf.\ \eqref{eq-symbol}) as pseudo-differential operator with symbol $q_+(s,x,\xi)$ then $L$ given by \eqref{fellergen} is on $C_c^\infty(\T \times\R^{d})$ a pseudo-differential operator and its symbol is given by
\begin{equation}\label{symol}
q_L(\tilde{x},\tilde{\xi}):= -i\sigma + q_+(s,x,\xi) 
\end{equation}
with $\tf{x}=(s,x)$ and $\tf{\xi} = (\sigma,\xi), s,\sigma \in \R,\ x,\xi \in \R^d.$ 
\end{Cor}
\begin{Proof}
By linearity the symbol of the sum of the operators is the sum of the symbols and for $s$ with covariable $\sigma$ the operator $\frac{d}{ds}$ corresponds to $-i\sigma$.
\end{Proof}

Note that for $\T=[0,\infty)$ one would have to take care of boundary terms.
Alternatively, in the setting of Corollary \ref{calsymbol}, a time inhomogeneous Markov process $X_t$ only defined for positive times, i.e.\ $\T=[0,\infty)$ can be extended onto $\T= \R$ by setting for $s<0$
\begin{equation}
q_+(s,x,\xi) := q_+(0,x,\xi).
\end{equation}

\begin{Rem}\label{exa-lrgen-spacetime}
One subtlety of Corollary \ref{calsymbol} is that it states that the generator of the time homogeneous process can be defined on $C_c^\infty(\T \times\R^{d})$, but it does not state that $C_c^\infty(\T \times\R^{d})$ is a subset of the domain of the generator of this process. To understand this, consider a Feller process on $\R$ whose sample path are deterministic with slope $\alpha>0$ below level $s_0$ and slope $\beta>0$ above level $s_0$, $\alpha\neq\beta$. Its generator is
\begin{equation} \label{disc-symbol}
L f(x) = l(x) f'(x) \quad \text { with } \quad l(x) := \begin{cases} \alpha &, x<s_0\\ \beta &, x\geq s_0 \end{cases}
\end{equation}
for all $f\in C^1_\infty(\R)$ with $f(s_0)=0$, where the last restriction is due to the requirement that $L: \dom(L) \to C_\infty.$ Thus $C_c^\infty(\R)$ is not a subset of the domain, nevertheless \eqref{disc-symbol} is well defined for $f \in C_c^\infty(\R)$ and $L$ has a representation as pseudo-differential operator with symbol $-i l(x) \xi.$ 

Finally note that the above discussion also applies to the process introduced in Example \ref{exa-lrgen} whose transformed process has the symbol
$$-i\sigma - il(s)\xi$$
where $(s,x)$ has covariable $(\sigma,\xi)$.
\end{Rem}

We close this section with a result which ensures that $C_c^\infty$ is in the domain of the transformed process.

\begin{Lem}\label{lem-dom} In the setting of Theorem \ref{calgenerator} let $C_c^\infty(\R^d) \subset \Dskript(A^+_s)$ for all $s$, and denote by $q_+(s,x,\xi)$ the symbol of $A^+_s$. If
$$s \mapsto q_+(s,x,\xi) \text{ is continuous for all }x,\xi$$
then 
$$C_c^\infty(\T\times\R^d) \subset \Dskript(L).$$
\end{Lem}
\begin{Proof}
Let $f\in C_c^\infty(\T\times\R^d)$ then $Lf$ is in $C_\infty(\T\times\R^d)$ by dominated convergence, and \eqref{fellergen} implies that $C_c^\infty(\T\times\R^d)$ is a subset of the domain of the pointwise generator. Finally the result follows, since for Feller semigroups the pointwise generator coincides with the generator (\cite[Lemma 31.7]{Sato99}, see also \cite[Lemma III.6.7\pages{, p.\ 241}]{Roge2000a}).
\end{Proof}

\section{Approximation of the process}
\label{secapprox}

Now we are going to show that a process $X_t$ with symbol $q_+$ given by \eqref{eq-symbol}  can -- under the assumptions of the following lemma  -- be approximated by Markov chains with time steps of size $\frac{1}{n}$. For each $n\in\N$ the approximating Markov chain $(Z^n(k))_{k\in\N_0}$ is defined by $Z^n(0):= X_0$ and transition kernels (from $x$ at time $k$ into $dy$ at time $k+1$) $\nu_{\frac{k}{n},x,\frac{1}{n}}(dy)$ where
\begin{equation}\label{approxtransition}
\int_{\R^d} e^{iy\xi} \ \nu_{s,x,\frac{1}{n}}(dy) = e^{ix\xi - \frac{1}{n} q_+(s,x,\xi)}.
\end{equation}

\begin{Lem} \label{conditions}
Let $(X_t)_{t\in\R}$ be a Markov process on $\R^d$ with corresponding Feller evolution system $U(s,t)$. Assume that $C_c^\infty(\R^d)$ is an operator core for the corresponding family of right generators $A_s^+$ , i.e.
$$\text{ for all $s$ the closure of } A_s^+\big|_{C_c^\infty} \text{ is } A_s^+.$$
Furthermore, assume that the symbol $-q_+(s,x,\xi)$ of $A_s^+\big|_{C_c^\infty}$ satisfies 
\begin{gather} \label{fellerevobound}
|q_+(s,x,\xi)| \leq c(1+|\xi|^2) \ \ \text{ for all } s,x,\xi,\\
\label{fellerevocont}
s\mapsto q_+(s,x,\xi) \ \ \text{ is continuous for all } x,\xi. 
\end{gather}
Then the test functions $C_c^\infty(\R^{d+1})$ are an operator core for the generator $L$ (given in Theorem \ref{calgenerator}) of the corresponding Feller process $\tf{X}_t$ and for the symbol $-\tf{q}(\tf{x},\tf{\xi})$ of $L\big|_{C_c^\infty(\R^{d+1})}$ exists a $c>0$ such that 
\begin{gather} \label{fellerbound}
|\tf{q}(\tf{x},\tf{\xi})|\leq c(1+|\tf{\xi}|^2) \ \ \text{ for all }\tf{x},\tf{\xi}.
\end{gather}
\end{Lem}
\begin{Proof} By Corollary \ref{calsymbol}
$$\tf{q}(\tf{x},\tf{\xi}) = -i\sigma + q_+(s,x,\xi)$$
holds and thus \eqref{fellerevobound} implies \eqref{fellerbound}. Furthermore, the core property follows by Lemma \ref{lem-dom} and linearity, since 
$$L = \frac{d}{ds} + A_s^+$$
and these operators act on different components.
\end{Proof}

Under the conditions of the Lemma \ref{conditions} the process $\tf{X}$ satisfies the assumptions of the approximation theorem in \cite{BoetSchi2009}. Thus, the process $\tf{X}_t$ is approximated by the Markov chains $(\tf{Y}^n([tk]))_{k\in\N_0}$ on $\R \times\R^d$ with $\tf{Y}^n(0):= (r,X_0)$ and transition kernel $\mu_{\tf{x},\frac{1}{n}}(d\tf{y})$ where
$$\int_{\R \times\R^d} e^{i\tf{y}\tf{\xi}} \ \mu_{\tf{x},\frac{1}{n}}(d\tf{y}) = e^{i\tf{x}\tf{\xi} - \frac{1}{n} \tf{q}(\tf{x},\tf{\xi})}.$$
Since the last $d$ coordinates of $\tf{X}_t$ started in $(0,x)$ coincide with $X_t$ started in $x$ we only need to check if we can simplify the above expression for these coordinates. 
Using $\tf{x} = (s,x),$ $\tf{y}= (r,y)$ and $\tf{\xi} = (\sigma,\xi)$ we get
\begin{eqnarray}
\tf{q}(\tf{x},\tf{\xi}) &= & - i\sigma + q(s,x,\xi),\\
e^{i\tf{x}\tf{\xi}-\frac{1}{n}\tf{q}(\tf{x},\tf{\xi})} &=& e^{i(s+\frac{1}{n})\tau+ix\xi - q(s,x,\xi)},\\
\mu_{\tf{x},\frac{1}{n}}(d\tf{y}) &=& \delta_{s+\frac{1}{n}}(dr)\times \nu_{s,x,\frac{1}{n}}(dy).
\end{eqnarray}
Thus the processes $\tf{X}$ can be approximated by the Markov chain defined by \eqref{approxtransition}. An approximations of this type is easily implemented for simulations, see \cite{BB2010} for an implementation of the time homogeneous case.

\bigskip

Finally we restate, using the transformation introduced in Section \ref{sectransfomation}, a result by Chernoff (\cite{Cher1974}, see also \cite[Theorem 2.5]{JacoPotr2004}) which shows that the approximation given above could also be used to construct Feller evolutions directly for a given family of probability measures $\nu_{s,x,\frac{1}{n}}$ as in \eqref{approxtransition}. The construction is formulated in terms of the operators 
$$V(s,h)g(x) := \int_{\R^d} g(y)\ \nu_{s,x,h}(dy)\ \ s\in\R, h\geq 0 \ \ \ g\in C_\infty(\R^d).$$
\begin{Cor} \label{approxconstruct}
If $(V(s,h))_{s\in\R,h\geq 0}$ is a family of strongly continuous linear contractions on $(C_\infty(\R^d),\|\cdot\|_\infty)$ which satisfies the following properties:
\begin{enumerate}[i)]
\item $V_{s,0} = id,$
\item the strong derivatives $\frac{d}{dh} V(s,h)\bigg|_{h=0} = V(s,0)'$ are densely defined,
\item 
$$\lim_{n\to\infty} \left\| V(s,\tfrac{t}{n})V(s+\tfrac{t}{n},\tfrac{t}{n})V(s+\tfrac{2t}{n},\tfrac{t}{n})\cdots V(s+\tfrac{(n-1)t}{n},\tfrac{t}{n}) g - U(s,t) g\right\|_\infty = 0$$
for all $g\in C_\infty(\R^d).$
\end{enumerate}
Then $(U(s,t))_{s,t\in\R, s\leq t}$ is a Feller evolution system on $C_\infty(\R^d)$ and its generators $A_s^+$ extend $V(s,0)'.$ Moreover the convergence in iii) is uniform for $t, s$ from compact intervals. 
\end{Cor}

\bibliographystyle{abbrv}

\end{document}